\documentclass[a4paper] {article}
[12pt]

\usepackage[top=3cm, bottom=3cm, left=3.8cm, right=3.8cm]{geometry}

\makeatletter
\renewcommand*\l@section{\@dottedtocline{1}{1.5em}{2.3em}}
\makeatother

\usepackage{amsfonts}
\usepackage{amssymb}
\usepackage[T1]{fontenc}

\usepackage{tikz}
\usetikzlibrary{calc}

\usepackage{CJK}
\usepackage{amsmath}

\usepackage{amsfonts}
\usepackage{amssymb}
\usepackage{amsthm}
\usepackage{amssymb}
\usepackage{enumerate}
\usepackage[calc]{picture}
\usepackage[all,cmtip]{xy}

\usepackage[mathscr]{eucal}
\usepackage{eqlist}

\usepackage{color}
\usepackage{abstract}
\usepackage[T1]{fontenc}

\theoremstyle{plain}
\newtheorem{theorem}{Theorem}
\newtheorem{proposition}[theorem]{Proposition}
\newtheorem{lemma}[theorem]{Lemma}

\newtheorem{example}[theorem]{Example}
\newtheorem{corollary}[theorem]{Corollary}

\theoremstyle{definition}
\newtheorem{definition}{Definition}

\usepackage{etoolbox}
\newtheoremstyle{myrem}
 {3pt}
 {3pt}
 {\normalsize}
 { }
 {\itshape}
 {:}
 { }
 {}

 \theoremstyle{myrem}
 \newtheorem{remark}{Remark}
 \appto\remark{\leftskip\parindent}
 \appto\remark{\rightskip\parindent}

\numberwithin{equation}{section}
\numberwithin{theorem}{section}

\begin{document}

\begin{center}
{\Large {\textbf {Differential  algebras   on   digraphs  and     parametrized     homology
}}}\footnotetext[1]{ This  is  an  improved  version  of  {\it  Differential algebras on digraphs and generalized path homology} (2021),
arXiv:2103.15870.  }
 \vspace{0.58cm}\\

Shiquan Ren,  Chong Wang 



\begin{quote}
\begin{abstract}
The theory of path homology for digraphs was developed by Alexander Grigor'yan,  Yong  Lin, Yuri Muranov,   and Shing-Tung Yau.  In this paper,  we  consider  the  differential  algebras  on digraphs and  define   the  parametrized  homology  of  digraphs  as  an analog  of the path homology.   We  prove  the  functoriality  of the  parametrized homology  of  digraphs  in  Theorem~\ref{th-88978}.   We  also  prove  some  K\"unneth-type  formulae  for  the  parametrized  homology   of  digraphs  in  Theorem~\ref{th-981328}.  
\end{abstract}
\end{quote}

\end{center}

\bigskip



\section{Introduction}

A  {\it digraph} $G$ is a couple $(V,D)$  where $V$  is a set and $D$ is a subset of $V\times V$ such that for any $(u,v)\in D$,  we have $u\neq v$.  
  During the 2010's,  Alexander Grigor'yan,  Yong  Lin, Yuri Muranov,   Shing-Tung Yau \cite{lin2,lin3,lin4,lin5},  
  Alexander Grigor'yan,  Yong  Lin,  Shing-Tung Yau \cite{lin1},     and  Alexander Grigor'yan,   Yuri Muranov,  Shing-Tung Yau  \cite{lin6}  developed  the theory of  path complexes and path homology for digraphs.

 Given a finite set $V$ and a non-negative integer $n\geq 0$,  an {\it $n$-path} on $V$  is a sequence $v_0v_1\ldots v_n$ of (not necessarily distinct)  $n+1$   elements in $V$.  Let $R$ be an integral domain such that $2$  has a multiplicative inverse $\frac{1}{2}$.  We let $\Lambda_n(V;R)$  be the free $R$-module generated by all the $n$-paths on $V$.  
For any element $v\in V$,  there  is an $R$-linear map   $\frac{\partial}{\partial v}$    from $\Lambda_n(V;R)$  to $\Lambda_{n-1}(V;R)$    whose explicit expression  is given by  (\ref{eq-1.1})   and  an  $R$-linear  map
  $dv$   from   from $\Lambda_n(V;R)$  to $\Lambda_{n+1}(V;R)$  whose   explicit  expression  is 
  given  by   (\ref{eq-2.2.1}).

Let  $R_*(V)$   be  the  free   $R$-module  generated  by   all  $\frac{\partial}{\partial  v}$  for  $v\in  V$.  
Let  $R^*(V)$   be  the  free   $R$-module  generated  by   all  $dv$  for  $v\in  V$.   By  some   algebraic   auxiliaries   on  differential algebras   and   parametrized  homology   (cf.   Section~\ref{s552}),   we  define  the  
parametrized  homology  of  digraphs  in   (\ref{def-aa23})  and  define  the  parametrized  cohomology  of  digraphs  in  (\ref{def-aa26}).

Let  $G$  be  a  digraph.    let  the  transitive-closure  $\Delta  G$  be  the  smallest  transitive  digraph containing  $G$  as  a  sub-digraph.  Then  $\Delta  G$  uniquely  exists. 
Let   the   transitive-interior  $\delta  G$  be  a  maximal  transitive  digraph  contained   in  $G$  as   a  sub-digraph.  
Note  that  $\delta  G$  exists while  may  not  be  unique.

In  Theorem~\ref{th-88978},  we  prove  the  functoriality  of  the  parametrized  homology  of  digraphs. 
As  a  consequence,  we  prove  the  functoriality  of   a    commutative  diagram  for  the  parametrized  homology
 of  $\delta  G$,  $G$,   and   $\Delta  G$,  in  Corollary~\ref{co-8912}.   
  In   Theorem~\ref{th-981328},  
we  prove  some  K\"unneth-type  formulae  for  the  parametrized  homology  of  joins  of  digraphs.  
As  a  consequence,  we  prove     a    commutative  diagram  for  the  K\"unneth-type  formulae  
 of  $\delta  G$,  $G$,   and   $\Delta  G$,  in  Corollary~\ref{co-uu76}.

  The  paper  is  organized  as  follows.   In   Section~\ref{s552},   we  give  some  algebraic  auxiliaries  on  
  differential  algebras  and  the   parametrized  homology. 
  In  Section~\ref{soo87},   we   define  the  transitive-closure  $\Delta  G$  and  the  transitive-interior   $\delta  G$  of  $G$.  In  Section~\ref{sss11qq},    we  study  the  differential  algebras  on  digraphs.  We  prove  Theorem~\ref{th-88978}   and   Theorem~\ref{th-981328}. 

\section{Differential algebras and their graded submodules}\label{s552}


\subsection{Parametrized differential algebras}\label{s1112}



Let $R$   be a commutative ring with unit such that    (1).  $2\in  R$  is  nonzero,  and  (2).  there exists  $1/2\in  R$ satisfying $2\cdot (1/2)=1$.
A  {\it $\mathbb{Z}$-graded  $R$-algebra}  is an  algebra $A$  over   $R$  such that there exists  a  family  $\{A_k\}_{k\in\mathbb{Z}}$ of      $R$-modules   with  both  of the followings   satisfied:
(1).  $A=\bigoplus_{k\in \mathbb{Z}}  A_k$  as $R$-modules,
(2).  $A_mA_n\subseteq A_{m+n}$  for any $m, n\in\mathbb{Z}$.
   An  {\it $\mathbb{N}$-graded    $R$-algebra}  is a $\mathbb{Z}$-graded $R$-algebra  such that  $A_k=0$ for any   integers  $k\leq 0$.
A {\it    differential  graded algebra (DGA)}  is  a   $\mathbb{N}$-graded  $R$-algebra  $A=\bigoplus_{k\geq 1}A_k$ with a sequence $d=\{d_k\}_{k\geq 0}$  of  morphisms of $R$-modules $d_k: A_k\longrightarrow  A_{k-1}$  (or  $d_k: A_k\longrightarrow  A_{k+1}$),  called the {\it differentials},  such  that   (1).  $d_{k}\circ d_{k+1}=0$ (or $d_{k+1}\circ d_{k}=0$) for  any  $k\geq 1$;  and  (2).   $d_{m+n}(e_me_n)=d_m(e_m)e_n+(-1)^{m}e_m d_n(e_n)$  for  any $m,n\geq 1$,  any  $e_m\in  A_m$  and  any  $e_n\in A_n$.
Let $X$  be an $R$-module.

\begin{definition}\label{def-881}
We  define  a {\it $X$-parametrized   differential  graded algebra ($X$-parametrized  DGA)}  to  be   a    $\mathbb{N}$-graded  $R$-algebra  $A=\bigoplus_{k\geq 1}A_k$  such that  (1).  for each $x\in X$
 there is a sequence $d(x)=\{d_k(x)\}_{k\geq 1}$ of  differentials $d_k(x): A_k\longrightarrow A_{k-1}$  (or  $d_k: A_k\longrightarrow  A_{k+1}$)  such that $A$  is a  DGA;  and  (2).  for any  $x,x'\in X$ and any  $k\geq  0$,   $d_k(x')d_{k+1}(x)=-d_k(x)d_{k+1}(x')$  (or  $d_{k+1}(x')d_{k}(x)=-d_{k+1}(x)d_{k}(x')$).  
 \end{definition}


Let  $\Lambda_R(X)=\bigoplus_{p\geq  0} \wedge^p X$  be  the  exterior algebra  over $R$  generated by  $X$.   Let  $A$   be  an  $X$-parametrized   DGA.    Then $\Lambda_R(X)$  acts  on  $A$
 naturally:
 \begin{eqnarray*}
 \Big(\sum r_{x_1\wedge  x_2\wedge \cdots \wedge  x_p}x_1\wedge  x_2\wedge \cdots \wedge  x_p\Big)(a)=\sum r_{x_1\wedge  x_2\wedge \cdots \wedge  x_p}d(x_1)\circ  d(x_2)\circ \cdots \circ  d(x_p) (a)
 \end{eqnarray*}
 for any  $a\in  A$.  In particular,  if $a\in A_k$  is an element  of homogeneous degree $k$,  where $k\geq 1$,  then  the righthand side of  the last equality  can  be written as
 \begin{eqnarray*}
 \sum r_{x_1\wedge  x_2\wedge \cdots \wedge  x_p}d_{k+1-p}(x_1)\circ  d_{k+2-p}(x_2)\circ \cdots \circ  d_k(x_p) (a).
 \end{eqnarray*}
  By  this  natural action of  $\Lambda_R(X)$  on $A$,  the  following  properties follow:
\begin{enumerate}[(1).]
\item
  For any $\alpha\in\wedge^p X$,  where $p$  is odd,  $\alpha\wedge \alpha=0$.  Thus   for  any $q=0,1,\ldots,p-1$ there is a  chain complex
\begin{eqnarray*}
\xymatrix{
 \cdots\ar[r]^-{d_{n+2}(\alpha,q)}
 & A_{q+(n+1)p} \ar[r]^-{d_{n+1}(\alpha,q)}
  &A_{q+np} \ar[r]^-{d_{n}(\alpha,q)}
  &A_{q+(n-1)p}
 \ar[r]^-{d_{n-1}(\alpha,q)}
 &\cdots \ar[r]^-{d_{1}(\alpha,q)}
  &A_q \ar[r]^-{d_{0}(\alpha,q)}
  &0
  }
\end{eqnarray*}
which will be denoted by  $C_*(A,\alpha,q)$.
Here if  $\alpha=\sum r_{x_1\wedge  x_2\wedge \cdots \wedge  x_p}x_1\wedge  x_2\wedge \cdots \wedge  x_p$ then  for each $n\geq  0$,
\begin{eqnarray*}
d_{n}(\alpha,q)= \sum r_{x_1\wedge  x_2\wedge \cdots \wedge  x_p}d_{q+(n-1)p+1}(x_1)\circ  d_{q+(n-1)p+2}(x_2)\circ \cdots \circ  d_{q+np}(x_p);
\end{eqnarray*}

\item
For any $\alpha\in\wedge^p X$,  where $p$  is odd,   and   any $\beta\in\wedge^r X$,  where $r$  is even,  $\alpha\wedge\beta=\beta\wedge \alpha$.  Thus  for any $q=0,1,\ldots,p-1$,  there is a  chain  map
\begin{eqnarray}\label{eq-2.3}
d_*(\beta):  C_*(A,\alpha,q)\longrightarrow  C_*(A,\alpha, q-r).
\end{eqnarray}
\end{enumerate}

\begin{definition}\label{def-98}
Let  $A$  be  an  $X$-parametrized   DGA.
For any $\alpha\in\wedge^p X$, where $p$  is odd,   any $q=0,1,\ldots,p-1$ and any positive  integer $n$,   we define the  $n$-th  {\it  $(\alpha,q)$-parametrized   homology} of  $A$   to be
\begin{eqnarray}\label{eq-88.a8}
H_n(A,\alpha,q)={\rm Ker}\big(d_n(\alpha,q)\big)/{\rm Im}\big(d_{n+1}(\alpha,q)\big).
\end{eqnarray}
  In addition,  for any $q=0,1,\ldots,p-1$,   we  define a  morphism of homology groups
\begin{eqnarray}\label{eq-88.b8}
\beta_*: H_n(A,\alpha,q)\longrightarrow H_n(A,\alpha,q-r),~~~ n=1,2\ldots
\end{eqnarray}
by applying  the homology functor to  (\ref{eq-2.3}).
\end{definition}


 Let  $A$  and  $A'$   be  two  $\mathbb{N}$-graded  $R$-algebras.  A  {\it morphism}  of $\mathbb{N}$-graded  $R$-algebras  from $A$  to $A'$,  written  $f:  A\longrightarrow A'$,    is a sequence
 $f=\{f_k\}_{k\geq  1}$   of   morphisms of  $R$-modules $f_k: A_k\longrightarrow A'_k$  such that $f_{m+n}(e_me_n)=f_m(e_m)f_n(e_n)$  for any  $m,n\geq 1$,  any  $e_m\in A_m$  and any  $e_n\in  A_n$.

 \begin{definition}\label{def-77}
 Let  $A$  and  $A'$   be  two  $X$-parametrized   DGAs.  For any  $x\in X$,  the differential of $A$  is $d(x)$  and the  differential of $A'$  is $d'(x)$.
  We  define  a   {\it morphism} of  $X$-parametrized   DGAs     from $A$  to $A'$   to  be  a   morphism  $f:  A\longrightarrow  A'$  of $\mathbb{N}$-graded  $R$-algebras  such  that  $d_k(x)\circ f_k=f_k \circ  d'_k(x)$ for any $x\in X$ and any $k\geq  1$.
 \end{definition}

 \begin{lemma}
 \label{th-1.1}
  For any $\alpha\in\wedge^p X$, where $p$  is odd,   any $\beta\in\wedge^r  X$, where $r$  is even,  and any $q=0,1,\ldots,p-1$,
both  the $(\alpha,q)$-parametrized  homology in  {\sc (\ref{eq-88.a8})}   and
 the morphism $\beta_*$  in  {\sc (\ref{eq-88.b8})}   are     functorial.
 \end{lemma}

 \begin{proof}
  Let $f: A\longrightarrow A'$  be a   morphism  of  $X$-parametrized   DGAs.   Then  $d'_n(\alpha)\circ  f_{np+q} =f_{(n-1)p+q} \circ  d_n(\alpha)$.
Thus  $f_{*p+q}$  is a  chain map from $C_*(A,\alpha,q)$ to $C_*(A',\alpha,q)$   which   induces a morphism from $H_*(A,\alpha,q)$  to $H_*(A',\alpha,q)$.   Hence the  $(\alpha,q)$-parametrized  homology  is  functorial.   Moreover,  $f_{np+q-r} \circ  d_n(\beta)= d'_n(\beta)\circ  f_{np+q}$.  
 Applying the  homology  functor,  we  obtain that
   $\beta_*$       is   functorial.
 \end{proof}

\subsection{Graded submodules of  parametrized   differential algebras and homology}\label{s1113}

Let  $X$  be an $R$-module.
Let  $A=\bigoplus_{k\geq 1}A_k$  be an  an $X$-parametrized   DGA.
 Let  $B=\bigoplus_{k\geq 1}B_k$  be  an $\mathbb{N}$-graded  sub-$R$-module  of $A$.  Then  for any $\alpha\in\wedge^p X$, where $p$  is odd,  any $\beta\in\wedge^r  X$, where $r$  is even,   and       any $q=0,1,\ldots,p-1$,
   there  is  a  chain  complex  ${\rm  Inf}_\bullet(B,A,\alpha,q)$  given by
 \begin{eqnarray*}
{\rm Inf}_n(B,A,\alpha,q)=B_{np+q}\cap \big(d_n(\alpha)\big)^{-1}(B_{(n-1)p+q}),~~~~~~n=1,2,\ldots
 \end{eqnarray*}
 and  a  chain  complex  ${\rm  Sup}_\bullet(B,A,\alpha,q)$  given by
 \begin{eqnarray*}
 {\rm  Sup}_n(B,A,\alpha,q)=B_{np+q}+ \big(d_{n+1}(\alpha)\big) (B_{(n+1)p+q}),~~~~~~n=1,2,\ldots
 \end{eqnarray*}
 Both ${\rm  Inf}_\bullet(B,A,\alpha,q)$  and ${\rm  Sup}_\bullet(B,A,\alpha,q)$   have  their   boundary  maps  $d_*(\alpha)$.
 By \cite[Proposition~2.4]{h1},  the  canonical inclusion $\iota: {\rm  Inf}_\bullet(B,A,\alpha,q)\longrightarrow {\rm  Sup}_\bullet(B,A,\alpha,q)$   is a quasi-isomorphism  of   chain  complexes.
 We define the {\it  $(\alpha,q)$-parametrized    homology groups}   of   $B$  by
 \begin{eqnarray}\label{eq-2.808}
 H_*(B,A,\alpha,q):=H_*\big({\rm  Inf}_\bullet(B,A,\alpha,q)\big)\overset{\iota}{\cong}  H_*\big({\rm  Sup}_\bullet(B,A,\alpha,q)\big).
 \end{eqnarray}
 By the commutative  diagram of chain  complexes
 \begin{eqnarray*}
 \xymatrix{
 {\rm Inf}_\bullet(B,A,\alpha,q)\ar[rr]^{d_r(\beta)}\ar[d]_{\iota}  &&{\rm Inf}_\bullet(B,A,\alpha,q-r)\ar[d]^{\iota}\\
 {\rm Sup}_\bullet(B,A,\alpha,q)\ar[rr]^{d_r(\beta)}  &&{\rm Sup}_\bullet(B,A,\alpha,q-r),
 }
 \end{eqnarray*}
   we  have an  induced  morphism of  homology  groups
 \begin{eqnarray}\label{eq-pqw}
  \beta_*: H_*(B,A,\alpha,q)\longrightarrow  H_*(B,A,\alpha,q-r).
 \end{eqnarray}

  \begin{lemma}
  \label{th-1.288}
  For any $\alpha\in\wedge^p X$, where $p$  is odd, any $\beta\in\wedge^r  X$, where $r$  is even,   and any $q=0,1,\ldots,p-1$,
 both
 the $(\alpha,q)$-parametrized    homology  in  {\sc (\ref{eq-2.808})}  and
 the morphism $\beta_*$   in {\sc (\ref{eq-pqw})}  are    functorial.
  \end{lemma}

\begin{proof}
 Let $f: A\longrightarrow A'$  be a   morphism  of  $X$-parametrized   DGAs.  Let  $B$  be  an $\mathbb{N}$-graded  sub-$R$-module  of $A$  and  let  $B'$  be  an $\mathbb{N}$-graded  sub-$R$-module  of $A'$  such  that $f(B)\subseteq  B'$.
Then $f: B\longrightarrow B'$  is a morphism of  $\mathbb{N}$-graded  $R$-modules.
By a diagram chasing,    
  $f_{* p+q}$  is a  chain map from $ {\rm  Sup}_\bullet(B,A,\alpha,q)$ to $ {\rm  Sup}_\bullet(B',A',\alpha,q)$.   This chain map  induces a morphism from $H_*(B,A,\alpha,q)$  to $H_*(B',A',\alpha,q)$.   It follows that  the $(\alpha,q)$-parametrized   path homology  is   functorial.
Morever,  we      have  a commutative  diagram  of  chain  complexes
 \begin{eqnarray*}
 \xymatrix{
 {\rm Sup}_\bullet(B,A,\alpha,q)\ar[rr] ^{d_*(\beta)} \ar[d]_{f_*}  &&   {\rm Sup}_\bullet(B,A,\alpha,q-r) \ar[d]^{f_*}\\
  {\rm Sup}_\bullet(B',A',\alpha,q)\ar[rr] ^{d_*(\beta)}  &&   {\rm Sup}_\bullet(B',A',\alpha,q-r).
 }
 \end{eqnarray*}
Applying the homology functor,  we  obtain that the morphism $\beta_*$  is  functorial.
\end{proof}

\begin{remark}
Let $B=A$  in Lemma~\ref{th-1.288}.  Then we obtain
 Lemma~\ref{th-1.1}.
\end{remark}

\subsection{K\"unneth-type  formulae for  parametrized   differential  algebras}\label{ss2.3}

Let $R$  be a principal ideal domain  with  unit $1$  such that       $2\in R$  is nonzero  and  there exists $1/2\in R$ satisfying $2\cdot (1/2)=1$.

\begin{theorem}{\sc \cite[Theorem~3B.5]{hatcher}}
\label{th-5.a1}
Let  $C$  and $C'$  be  chain  complexes  of  free $R$-modules.
For any  integer $n$,  there  is  a  natural  short exact sequence
 \begin{eqnarray}
&0\longrightarrow  \bigoplus_{t+s=n} H_t(C)\otimes H_s(C')\longrightarrow H_n(C\otimes C')\nonumber\\
& \longrightarrow \bigoplus_{t+s=n }{\rm Tor}_R(H_t(C),H_{s-1}(C'))\longrightarrow 0.
 \label{eq-5.a987}
 \end{eqnarray}
\end{theorem}

Let  $X$  be an $R$-module.
Let  $A=\bigoplus_{k\geq 1}A_k$  and   $A'=\bigoplus_{k\geq 1}A'_k$   be two   $X$-parametrized   DGAs  of  free $R$-modules.  Then  their tensor product
\begin{eqnarray*}
A\otimes A'=\bigoplus_{k\geq 1} \Big(\bigoplus_{t+s=k\atop t,s\geq 1}A_t\otimes A'_s \Big)
\end{eqnarray*}
is a graded  free $R$-module  given by  $(A\otimes A')_k= \bigoplus_{t+s=k\atop t,s\geq 1} A_t\otimes A'_s$  for  each  $k\geq 1$.
 For any $x\in X$,  define  the differential
 \begin{eqnarray}\label{eq-tpd}
 d(x): A_t\otimes A'_s\longrightarrow  A_{t-1}\otimes A'_s + A_t\otimes A'_{s-1}
 \end{eqnarray}
  by the Leibniz rule
  \begin{eqnarray}\label{eq-lnr}
 d(x)(c\otimes c')=d(x)(c)\otimes c' + (-1)^t c\otimes d(x)(c'),
 \end{eqnarray}
 where  $c\in A_t$  and  $c'\in A'_s$.
   Then
 (\ref{eq-tpd})  gives a differential
 \begin{eqnarray*}
 d(x): (A\otimes A')_k\longrightarrow (A\otimes A')_{k-1},  ~~~ k=1,2,\ldots
 \end{eqnarray*}
By  Definition~\ref{def-881},  $A\otimes A'$  is   an $X$-parametrized  DGA 
 of  free $R$-modules.

\begin{lemma}\label{le-kf}
For any $x\in X$  and  any  nonnegative  integer $n$,  there  is  a  natural  short exact sequence
   \begin{eqnarray}
 &0\longrightarrow  \bigoplus_{t+s=n} H_t(A,x,0)\otimes H_s(A',x,0)\longrightarrow H_n(A\otimes A', x,0)
 \nonumber\\
 &\longrightarrow \bigoplus_{t+s=n }{\rm Tor}_R(H_t(A,x,0),H_{s-1}(A',x,0))\longrightarrow 0.
 \label{eq-5.a98898}
 \end{eqnarray}
\end{lemma}

\begin{proof}
In Theorem~\ref{th-5.a1},  let $C$  be the chain complex $A$  and
 let $C'$  be the chain complex $A'$.
 Here  the boundary map  of $A$  is  $d_n(x): A_n\longrightarrow A_{n-1}$, $n=1,2,\ldots$,    and
  the boundary map  of $A'$  is  $d_n(x): A'_n\longrightarrow A'_{n-1}$, $n=1,2,\ldots$. 
 Then  $C\otimes C'$ is the chain complex $A\otimes A'$
 with the boundary map    $d_n(x): (A\otimes A')_n\longrightarrow  (A\otimes A')_{n-1}$  given by the Leibniz rule  (\ref{eq-lnr}).
 The lemma follows from Theorem~\ref{th-5.a1}.
\end{proof}

The next corollary  follows from Lemma~\ref{le-kf}.

\begin{corollary}\label{le-kf11}
Let $R$  be a field $\mathbb{F}$  such that  ${\rm char}(\mathbb{F})\neq 2$.   Then for any $x\in X$,  we have an isomorphism of graded $\mathbb{F}$-vector spaces $H_*(A\otimes A',x,0)\cong  H_*(A,x,0)\otimes H_*(A',x,0)$.  \qed
\end{corollary}

Let $x_1, x_2\in X$.  The  chain complex
 $C_*(A,x_1,0)$  is  given by
 \begin{eqnarray*}
  \xymatrix{
  \cdots \ar[r]^-{d_{n+1}(x_1)}
  & A_n  \ar[r]^-{d_{n }(x_1)}
  &A_{n-1} \ar[r]^-{d_{n-1}(x_1)}
  &\cdots \ar[r]^-{d_{1}(x_1)}
  &A_0.   }
  \end{eqnarray*}
  By Definition~\ref{def-881}~(2),  for each  nonnegative integer $n$ the homomorphism
  $d_n(x_2)$ sends
  ${\rm Im}(d_{n+1}(x_1))$ to  ${\rm Im}(d_{n}(x_1))$
  and sends  ${\rm Ker}(d_{n}(x_1))$ to  ${\rm Ker}(d_{n-1}(x_1))$.
  Thus $d_n(x_2)$ gives a homomorphism  of $R$-modules
\begin{eqnarray*}
   d_n(x_2): H_n(A,x_1,0)\longrightarrow H_{n-1}(A,x_1,0).
\end{eqnarray*}
   This gives  a chain complex $C_*\big(H(A,x_1,0),x_2,0\big)$
    by
  \begin{eqnarray*}
   \xymatrix{
 \cdots \ar[r]^-{d_{n+1}(x_2)}
  & H_n(A,x_1,0)  \ar[r]^-{d_{n }(x_2)}
  &H_{n-1}(A,x_1,0)\ar[r]^-{d_{n-1}(x_2)}
  &\cdots \ar[r]^-{d_{1}(x_2)}
  &H_0(A,x_1,0).
  }
 \end{eqnarray*}
 Denote the homology groups of $C_*\big(H(A,x_1,0),x_2,0\big)$  by
 $H_* (A,x_1, x_2,0 )$.

Generally,  let $m\geq  2$  be an     integer.  Let $x_1,x_2,\ldots,x_m\in X$.
 Inductively   we have a sequence of chain complexes
 \begin{eqnarray*}
 C_*\big(H(A,x_1,x_2,\ldots,x_{i-1},0),x_i,0\big), ~~~ i=1,2,\ldots,m
 \end{eqnarray*}
 whose homology groups are denoted by
 \begin{eqnarray}\label{eq-comp88}
 H_*(A, x_1, x_2,\ldots,x_i, 0 ), ~~~ i=1,2, \ldots,m.
 \end{eqnarray}
 We call (\ref{eq-comp88}) the {\it $(x_1, x_2,\ldots,x_i)$-composed  parametrized  homology} of $A$.

\begin{corollary}\label{co-composed}
Let $R$  be a field $\mathbb{F}$  such that   ${\rm char}(\mathbb{F})\neq 2$.   Then for any
positive integer $m$ and any $x_1,x_2,\ldots,x_m\in X$,
we have an isomorphism of graded $\mathbb{F}$-vector spaces
\begin{eqnarray}\label{eq-kn988}
H_*(A\otimes A',x_1, x_2,\ldots,x_m,0)\cong  H_*(A,x_1, x_2,\ldots,x_m,0)\otimes H_*(A',x_1, x_2,\ldots,x_m,0).
\end{eqnarray}
\end{corollary}
\begin{proof}
We prove by an induction on $m$.  When $m=2$,  (\ref{eq-kn988}) follows by Corollary~\ref{le-kf11}.   Suppose (\ref{eq-kn988})  holds  for $m-1$.  Then
\begin{eqnarray*}
&&H_*(A\otimes A',x_1, x_2,\ldots,x_m,0)\\
&=& H\Big(C_*\big(H(A\otimes A',x_1,x_2,\ldots,x_{m-1},0),x_m,0\big)\Big)\\
&\cong &H\Big(C_*\big(H(A,x_1,x_2,\ldots,x_{m-1},0)\otimes H(A',x_1,x_2,\ldots,x_{m-1},0) ,x_m,0\big)\Big)\\
&=&H\Big(C_*\big(H(A,x_1,x_2,\ldots,x_{m-1},0),x_m,0\big)\\
&&
\otimes C_*\big(H(A,x_1,x_2,\ldots,x_{m-1},0) ,x_m,0\big)\Big)\\
&\cong &H\Big(C_*\big(H(A,x_1,x_2,\ldots,x_{m-1},0),
x_m,0\big)\Big)\\
&&\otimes H\Big(C_*\big(H(A,x_1,x_2,\ldots,x_{m-1},0) ,x_m,0\big)\Big)\\
&=&H_*(A,x_1, x_2,\ldots,x_m,0)\otimes H_*(A',x_1, x_2,\ldots,x_m,0).
\end{eqnarray*}
Thus we obtain (\ref{eq-kn988})    for any  positive integer $m$.
\end{proof}

\subsection{ K\"unneth-type  formulae for   graded submodules of
parametrized  differential  algebras}\label{ss2.4}

Let $R$  be a principal ideal domain  with  unit $1$.

\begin{theorem}{\sc \cite[Theorem~1.1]{wrl}}
\label{th-5.a22}
Let  $C=\bigoplus_{i\geq 0} C_i$  and $C'=\bigoplus_{i\geq 0}C'_i$  be  two  chain  complexes  of  finitely-generated free $R$-modules.
For each  nonnegative integer $i$  we  let $b_i$  be the set of the generators of $C_i$ and  let  $b'_i$  be the set of the generators of $C'_i$.
Let $e_i$  be a subset of $b_i$ and let  $e'_i$  be a subset of $b'_i$.
Let  $D=\bigoplus_{\geq 0} D_i$  and $D'=\bigoplus_{i\geq 0} D'_i$  be  the graded free $R$-modules   where each $D_i$  is generated by $e_i$ and each $D'_i$  is generated by $e'_i$.
Then for any nonnegative  integer $n$,  there  is  a  natural  short exact sequence
 \begin{eqnarray}
&0\longrightarrow  \bigoplus_{t+s=n} H_t(D,C)\otimes H_s(D',C')\longrightarrow H_n(D\otimes D',C\otimes C')
 \nonumber\\
 &\longrightarrow \bigoplus_{t+s=n }{\rm Tor}_R(H_t(D,C),H_{s-1}(D',C'))\longrightarrow 0.
 \label{eq-5.a98788}
 \end{eqnarray}
 Here $H_*(D,C)$  is the  path homology of $D$   i.e. the homology groups of the  infimum chain complex ${\rm Inf}_*(D,C)$,  $H_*(D',C')$  is the  path homology of $D'$   i.e. the homology groups of the  infimum chain complex ${\rm Inf}_*(D',C')$,  and  $H_*(D\otimes D',C\otimes C')$  is the  path homology of $D\otimes D'$   i.e. the homology groups of the  infimum chain complex ${\rm Inf}_*(D\otimes D',C\otimes C')$.
\end{theorem}

  Suppose       $2\in R$  is nonzero  and  there exists $1/2\in R$ satisfying $2\cdot (1/2)=1$.
Let  $X$  be an $R$-module.
Let  $A=\bigoplus_{k\geq 1}A_k$  and  $A'=\bigoplus_{k\geq 1}A'_k$
be  two  $X$-parametrized   DGAs.
For each positive integer $k$,
suppose  $A_k$  and  $A'_k$  are  finitely-generated free $R$-modules.
  Let $a_k$  be the   set of  the  generators of $A_k$  and  let   $a'_k$  be the  set of the generators of $A'_k$.
    Let $b_k$  be a subset of $a_k$  and let  $b'_k$  be a subset of $a'_k$.
 Let $B_k$  be  the  finitely-generated  free  $R$-module generated by $b_k$ and
 let  $B'_k$  be  the  finitely-generated free  $R$-module generated by $b'_k$.
Then $B=\bigoplus_{k\geq 1} B_k$  is  an $\mathbb{N}$-graded  sub-$R$-module  of $A$  and $B'=\bigoplus_{k\geq 1} B'_k$  is  an $\mathbb{N}$-graded  sub-$R$-module  of $A'$.

\begin{lemma}\label{le-ceh2}
 For any $x \in X$  and  any    nonnegative  integer $n$,   there  is  a  natural  short exact sequence
 \begin{eqnarray}
&0\longrightarrow  \bigoplus_{t+s=n} H_t(B,A, x,0)\otimes H_s(B',A', x,0)\longrightarrow H_n(B\otimes B',A\otimes A', x,0)
 \nonumber\\
 &\longrightarrow \bigoplus_{t+s=n }{\rm Tor}_R(H_t(B,A, x,0),H_{s-1}(B',A', x,0))\longrightarrow 0.
 \label{eq-5.a9998}
 \end{eqnarray}
\end{lemma}

\begin{proof}
In  Theorem~\ref{th-5.a22},
let $C$  be  the chain complex $C_*(A, x,0)$,
let  $C$  be  the chain complex $C_*(A', x,0)$,
 let  $D$  be  the  graded $R$-module  $B$
 and let  $D'$  be  the graded $R$-module $B'$.
  Then  (\ref{eq-5.a9998}) follows from  (\ref{eq-5.a98788}).
\end{proof}

\section{Digraphs  and     transitive  digraphs}\label{soo87}

A  digraph  $G$  consists  of  a   set  $V$   of  vertices,   a  set  $D$  of  directed  edges,   and  maps  
$s,t:  D\longrightarrow   V$   assigning  to  each  $e\in  D$  its  source  vertex  $s(e)$  and  its  target  
vertex  $t(e)$.    We  assume  $s(e)\neq  t(e)$  for  any  $e\in  D$.  
 Throughout  this  paper,   we  assume  that  there   is   no  multiple  edges  in  $G$,  that  is,   for  any   $e,e'\in  D$,   if   $s(e)=s(e')$  and  $t(e)=t(e')$,   
 then   $e=e'$.

  A   {\it  morphism}  of  digraphs    $p=(p_0,p_1)$   from   $G=(V,D,s,t)$  to  $G'=(V',D',s',t')$  consists  of  two  maps 
 $p_0:   V\longrightarrow  V'$  and   $p_1:  D\longrightarrow  D'$  such  that  $s'\circ  p_1=  p_0\circ  s$  and 
 $t'\circ  p_1= p_0\circ  t$.  
 We  say  that   $G$  is  a  sub-digraph  of  $G'$  if  there  exists  a  morphism  of  digraphs  $p=(p_0,p_1)$  such  that  both  $p_0$  and  $p_1$  are  injective.

We  call  a  digraph  $G=(V,D,s,t)$   {\it  transitive}  if  for  any   $e_1,e_2\in  D$  such  that  $t(e_1)=s(e_2)$,  there exists  
$e\in  D$  such  that   $s(e)=s(e_1)$  and   $t(e)=t(e_2)$.

Let   $G=(V,D,s,t)$   be  a  fixed  digraph.   Let  $\Delta  G=(V,\Delta  D, \Delta  s, \Delta  t)$   be  the  minimal   transitive  digraph   such  that  
$G$  is  a  sub-digraph  of   $\Delta  G$.   That  is,
\begin{eqnarray}\label{eq-2.18a}
\Delta  G=\bigcap _{G{\rm~is~a~sub{\text-}digraph~of~}G' \atop G'  {\rm  ~is~transitive}} G'. 
\end{eqnarray}
Note that  the intersection  of  certain  transitive digraphs  is  still  a  transitive  digraph.  Thus  (\ref{eq-2.18a})  makes sense.     
  Let  $\delta  G=(V,\delta  D, \delta  s, \delta  t)$  be  a  maximal   digraph  such  that   $\delta  G$  is  
a  sub-digraph  of   $G$.   We  call  $\Delta  G$  the  {\it   transitive-closure  digraph}  of  $G$  and  call  $\delta  G$  the  {\it     transitive-interior  digraph}  of  $G$.    Note that  $\Delta  G$  is  unique  while  $\delta  G$  may  not  be  unique.   Choose an  arbitrary  transitive-interior digraph  $\delta  G$  of  $G$.  
We  have  $\delta  D\subseteq  D\subseteq  \Delta  D$.   Moreover,  
$\delta   s$  is  the restriction  of  $s$  to   $\delta  D$  and   $s$  is  the  restriction  of  $\Delta  s$  to  $D$.

  Let  $|G|=(V,E)$   be  the  graph  where   $E$  
is  the  collection of  the  unordered  pairs  $uv$   such  that  there exists  $e\in  D$  with 
$s(e)=u$  and  $t(e)=v$  or  $s(e)=v$  and  $t(e)=u$.   We  call  $|G|$  
the  {\it  underlying graph}   of  $G$.  
 The  digraph  $G$  
can  be  obtained  from    $|G|$  by  assigning  a  direction  or  both  directions to each edge  in  $E$.    
   We  call   a  sub-digraph 
$G_0$  a  {\it connected  component}  of  $G$  if      $|G_0|$    is   a   connected  component  of   $|G|$.  

Let  $G$  be  a  digraph.   Let  $G=\bigsqcup _i  G_i$  be  its  decomposition into  connected  components.  
Then  $\Delta  G=\bigsqcup _i  \Delta(G_i)$  is  the  decomposition  of  $\Delta  G$  into  connected  components,  where  $\Delta (G_i)$  is  the  transitive-closure  of  $G_i$.   Moreover,  for  any  transitive-interior  digraph  $\delta  G$  of  $G$,  we  have  a  decomposition $\delta  G=\bigsqcup _i  \delta(G_i)$,    where  each  $\delta (G_i)$  is  
a  transitive-interior  of  $G_i$.  

 Let   $p:  G\longrightarrow   G'$  be  a  morphism  of   digraphs.   For  any  transitive-interior  $\delta  G$  of  $G$,  since  its  image  under   $p$   must  be  transitive,  there  exists  a  transitive-interior  $\delta  G'$  of  $G'$  such  that  the  following      diagram   commutes
\begin{eqnarray*}
\xymatrix{
\delta  G\ar[r]^-{\delta  p} \ar[d]  & \delta  G' \ar[d]\\
G\ar[r]^-p \ar[d] &G'\ar[d]\\
\Delta  G\ar[r]^-{\Delta  p}&  \Delta  G'.    
}
\end{eqnarray*}
Here  $\delta  p$  is  the  restriction  of  $p$  to  $\delta  G$  and  $p$  is  the  restriction  of   $\Delta  p$  
to  $G$.   Consider  the  decompositions    $\bigsqcup _i  G_i$   of   $G$ and   $\bigsqcup _j  G'_j$   of   $G'$  into  connected  components.   Then   for  each  component  $G_i$  there  exists  a  component  $G'_j$  such  that  
the  following      diagram   commutes
\begin{eqnarray*}
\xymatrix{
\delta  (G_i)\ar[r]^-{\delta  p} \ar[d]  & \delta  (G'_j) \ar[d]\\
G_i\ar[r]^-p \ar[d] &G'_j\ar[d]\\
\Delta  (G_i)\ar[r]^-{\Delta  p}&  \Delta  (G'_j).    
}
\end{eqnarray*}

\begin{example}\label{ex-2.1a}
Let  $V=\{v_0,v_1,v_2, v_3\}$.   Let  $D=\{e_1,e_2,e_3,e_4,e_5\}$.   Let  
\begin{eqnarray*}
&s(e_1)= v_0, ~~~~~~ t(e_1)=v_1, \\
&s(e_2)= v_1, ~~~~~~ t(e_2)=v_2, \\
&s(e_3)=v_2, ~~~~~~ t(e_3)=v_0, \\
&s(e_4)=v_2, ~~~~~~ t(e_4)= v_1,\\
&s(e_5)=v_0, ~~~~~~ t(e_5)=v_3. 
\end{eqnarray*}
Then 
\begin{eqnarray*}
 \Delta  D &=& \{e_1,e_2,e_3,e_4,e_5, v_0 v_2,  v_1v_0, v_1v_3, v_2v_3  \}, \\
 \delta  D &=& \{e_1,e_2,e_3,e_4, v_0 v_2,  v_1v_0\}
\end{eqnarray*}
and
\begin{eqnarray*}
&(\Delta s)(v_0v_2)=v_0,  ~~~~~~  (\Delta  t)(v_0v_2)=v_2, \\
&(\Delta s)(v_1v_0)=v_1,  ~~~~~~  (\Delta  t)(v_1v_0)=v_0, \\
&(\Delta s)(v_1v_3)=v_1,   ~~~~~~  (\Delta t)(v_1,v_3)=v_3, \\
&(\Delta s)(v_2v_3)=v_2,   ~~~~~~  (\Delta t)(v_1,v_3)=v_3. 
\end{eqnarray*}
\end{example}

\begin{example}\label{ex-2.2a}
Let  $G=(V,D,s,t)$  be  given  in  Example~\ref{ex-2.1a}.  
\begin{enumerate}[(i).]
\item
Let   $G_\infty=(V_\infty,D_\infty,s_\infty,t_\infty)$  be  the  digraph  given  by 
$V_\infty=\mathbb{Z}$,  $D_\infty=\{(k,k+1)\mid  k\in\mathbb{Z}\}$,   
$s_\infty(k,k+1)=k$  and   $t_\infty(k,k+1)=k+1$.  
 Let  $p:  G_\infty\longrightarrow  G$  be  a  morphism  of  digraphs  given  by 
\begin{eqnarray*}
  p_0(3l+r)= r   ~~~{\rm   and } ~~~
   p_1(3l+r-1,  3l+ r)=  e_r  
\end{eqnarray*}
where  $l\in\mathbb{Z}$  and   $r=0,1,2$.  
We  have  
\begin{eqnarray*}
&  \Delta  (D_\infty)=  \{(a,b)\mid   a,b\in\mathbb{Z},   a<b\},\\
&  s(a,b)=a,~~~~~~    t(a,b)=b,\\
&  \Delta  p_1 (a,b)=   (p_0(a), p_0(b)), 
\end{eqnarray*}
where  we  use  the  notation  $e_1=v_0v_1$,  $e_2=v_1v_2$,  $e_3=v_2v_0$,   $e_5=v_2v_1$. 
Moreover,  
\begin{eqnarray*}
&  \delta  (D_\infty)=  \{(2k,2k+1)\mid   k\in\mathbb{Z}\},  \\
&  \Delta  p_1 (2k,2k+1)=   (p_0(2k), p_0(2k+1))  
\end{eqnarray*}
or
\begin{eqnarray*}
&  \delta  (D_\infty)=  \{(2k-1,2k)\mid   k\in\mathbb{Z}\},\\
&  \Delta  p_1 (2k-1,2k )=   (p_0(2k-1), p_0(2k)).    
\end{eqnarray*}
\item
Let  $m$  be  a  positive  integer.  
Let   $G_m=(V_m,D_m,s_m,t_m)$  be  the  digraph  given  by 
$V_m=\{k\in\mathbb{Z}\mid   0\leq  k\leq  3m-1\}$,  $D_\infty=\{(k,k+1)\mid  k\in\mathbb{Z}, 0\leq k\leq  3m-2\}\cup\{(3m-1,0)\}$.  
 Let  $p:  G_m\longrightarrow  G$  be  a  morphism  of  digraphs  given  by 
\begin{eqnarray*}
 & p_0(3l+r)= v_r,\\
  &  p_1(3l+r-1,  3l+ r)=  e_r,  ~~~p_1(3m-1,0)= e_3    
\end{eqnarray*}
where  $l=0,1,\ldots, m-1$  and   $r=0,1,2$.  
 We  have  
\begin{eqnarray*}
&  \Delta  (D_m)=  \{(a,b)\mid   a,b\in  V_m\},\\
&  s(a,b)=a,~~~~~~    t(a,b)=b,\\
&  \Delta  p_1 (a,b)=   (p_0(a), p_0(b)).  
\end{eqnarray*}
Moreover,  
\begin{eqnarray*}
  \delta  (D_m)= 
\begin{cases}
 \{(0,1),  (2,3),\ldots, (3m-1,0) \}   &{\rm~if~} m{\rm~is ~odd}, \\
 \{(0,1), (2,3), \ldots, (3m-2,3m-1)\}  &{\rm~if~} m{\rm~is ~even}
 \end{cases}
\end{eqnarray*}
or
\begin{eqnarray*}
  \delta  (D_m)= 
\begin{cases}
 \{(1,2),  (3,4),\ldots, (3m-2,3m-1) \}   &{\rm~if~} m{\rm~is ~odd}, \\
 \{(1,2), (3,4), \ldots, (3m-1,0)\}  &{\rm~if~} m{\rm~is ~even}. 
 \end{cases}
\end{eqnarray*}
The  morphism of digraphs  $\delta  p$  is  the restriction  of  $p$ to  $\delta  (G_m)$. 
\end{enumerate}
\end{example}
 
We  say  that  $G=(V,D,s,t)$  is  a  {\it  complete  digraph}   if   for  any  ordered  pair    $v,u\in  V$  with  $u\neq  v$,  there   (uniquely)  exists     $e\in  D$  such  that    $s(e)=v$  and  $t(e)=u$.    
It  is  direct  to  see  that   $G$  is  a  complete  digraph  iff  $|G|$   is  a  complete  graph  and  each  edge  of  $G$  is  assigned  with  both  directions.    
Let  $G$  be a  complete  digraph.  Then   $\delta  G=\Delta  G=G$.  In this case,  $\delta  G$  is  unique.   Let  $p:  G\longrightarrow  G'$    be  a  morphism  of  digraphs  
where  both  $G$  and  $G'$  are  complete  digraphs.   Then   $\delta  p = \Delta  p  =p$.  

   \begin{example}
Let  $G$  be  a  digraph.  A  {\it  cycle}  of  $G$  is  a  sequence $e_1e_2\ldots  e_n e_0$,  
   where  $ e_0,  e_1, \ldots,  e_n\in  D$,   such  that   
$s(e_k)=t(e_{k-1})$,  $k=1,\ldots, n$  and  $s(e_0)=t(e_n)$.       It  is  direct to prove  that  if  there  is  a  cycle  of  $G$  such  that  each   $v\in  V$  is  either  the  source or  the   tail  of  certain  $e_k$,  then  $\Delta  G$  is  a  complete  digraph  on   $V$.        
 \end{example}

We  say  that  $G=(V,D,s,t)$  is  a  {\it  semi-complete  digraph}  if   for  any  ordered  pair    $v,u\in  V$  with  $u\neq  v$,   (1).  either   (1.a). there   (uniquely)   exists    $e\in  D$  such  that      $s(e)=v$  and  $t(e)=u$  or  (1.b).  there   (uniquely)   exists  $e\in  D$  such  that   $s(e)=u$  and  $t(e)=v$;  and  (2).   (1.a)  and  (1.b)  cannot  both  be  true.   It  is  direct 
to  see   that   $G$  is  a  semi-complete  digraph  iff  $|G|$  is  a  complete  graph  and  each  edge  of   $G$   
 is  assigned  with  exactly  one  direction.  

 \begin{example}
Consider  the  line  digraph  $I_n=([n],  D,s,t)$,  where  $[n]=\{0,1,\ldots,n\}$,    $D=\{(k-1,k)\mid k=1, \ldots, n \}$,
 $s(k-1,k)=k-1$  and   $t(k-1,k)=k$.       Then    
$\Delta  (I_n)$   is   a   semi-complete   digraph  on   $[n]$  whose  set of  directed  edges  are the  collection  of  
 $(a,b)$,  $0\leq  a<b \leq  n$.       
  \end{example}

\section{The  parametrized   homology of digraphs}\label{sss11qq}

 Let  $G=(V,D,s,t)$  be  a  digraph. 
 We  say  that  $G$  is  {\it  locally  finite}  if  for  any  $v\in  V$,  there exists  finitely  many  $e\in  D$  such  that 
 either  $s(e)=v$  or  $t(e)=v$.   
 A  $0$-path  on  $G$  is  a  vertex  of  $G$.   
 Let  $n$  be  a  positive   integer.  
 An  {\it  elementary  $n$-path} on $G$       is a sequence $v_0v_1\ldots v_n$ of (not necessarily distinct)  $n+1$   vertices  in $V$  such  that     for  each  $1\leq  i \leq  n$,  either   $v_{i-1}=v_i$  or  
  there  (uniquely)  exists  $e\in  D$  satisfying  $s(e)=v_{i-1}$  and  $t(e)=v_i$.   
 Let  $\Lambda_n(G;R)$  be  the free $R$-module  generated  by  all the   elementary  $n$-paths  on  $G$. 
 Let   $\Lambda_*(G;R)= \bigoplus_{n\geq   0}  \Lambda_n(G;R)$.

 Let  $\Delta[V]$   be  the  complete  digraph  on  $V$.   
 An  elementary $n$-path  on  $\Delta[V]$  is  an  arbitrary  sequence  $v_0v_1\ldots  v_n$  such  that  $v_{i-1}\neq  v_i$  for  each
 $1\leq  i\leq  n$.     
For  any  nonnegative  integers  $n$  and   $m$,  the  join   is  an   $R$-bilinear  map 
$
*:  \Lambda_n(V;R)\times  \Lambda_m(V;R)\longrightarrow  \Lambda_{m+m+1}(V;R)
$
  given  by   
\begin{eqnarray*}
v_0\ldots  v_n  *  u_0\ldots  u_m  =   v_0\ldots  v_n  u_0\ldots u_m.  
\end{eqnarray*}
For  each  $n\geq  0$,  define  an  $R$-linear  map  $d_n:  \Lambda_n(V;R)\longrightarrow  \Lambda_{n-1}(V;R)$  by 
\begin{eqnarray*}
d_n  (v_0v_1\ldots  v_n) =\sum_{i=0}^n  (-1)^i  v_0\ldots \widehat {v_i}\ldots  v_n.   
\end{eqnarray*}
Here  we  write  $  \Lambda_{-1}(V;R)= 0$.    
 It  is  direct to verify that  $d_n\circ  d_{n+1}=0$   and   $d(v_0\ldots  v_n  *  u_0\ldots  u_m)=d_n (v_0\ldots v_n)  *  u_0\ldots u_m + (-1)^{n+1} v_0\ldots  v_n * d_m (u_0\ldots  u_m)$  for  any  nonnegative  integers  $m$  and  $n$.    Thus  
 $\Lambda_*(V;R)= \bigoplus_{n\geq   0}  \Lambda_n(V;R)$   is   a   DGA (in  Section~\ref{s1112},  let  $A_k=\Lambda_{k-1}(V;R)$).

 Let  $R(V)$  be  the  free  $R$-module  generated  by  all the  vertices  in  $V$.  
We  generalize  the  DGA  structure  in  the  last  paragraph  and  define    $R(V)$-parametrized  DGA  structures  on  $\Lambda_*(V;R)$.   
 For any $v\in V$,  we define   a graded  $R$-linear map
\begin{eqnarray*}
\frac{\partial}{\partial v}: ~~ \Lambda_n(V;R)\longrightarrow \Lambda_{n-1}(V;R),~~~ n\geq 0,
\end{eqnarray*}
by letting
\begin{eqnarray}\label{eq-1.1}
\frac{\partial}{\partial v} \Big(v_0v_1 \ldots v_n\Big)=\sum_{i=0}^n (-1)^i  \delta(v,v_i) v_0\ldots \widehat{v_i}\ldots v_n.
\end{eqnarray}
Here in (\ref{eq-1.1}),   for  any  vertices  $u,v\in V$  we use the notation  $\delta(u,v)=1$  if $u=v$  and  $\delta(u,v)=0$  if  $u\neq v$.   Particularly,  if   $v_0,v_1,\ldots,v_n$ are distinct vertices in $V$,  then   we have the followings: 
\begin{itemize}
\item
if $v_i=v$ for some $0\leq i\leq n$, then 
\begin{eqnarray*}
\frac{\partial}{\partial v} \Big(v_0 v_1\ldots v_n\Big)=(-1)^i   v_0\ldots \widehat{v_i}\ldots v_n; 
\end{eqnarray*}
 \item
 if  $v_i\neq v$  for any $0\leq i\leq n$, then 
\begin{eqnarray*}
\frac{\partial}{\partial v} \Big(v_0v_1\ldots v_n\Big)=0.  
\end{eqnarray*}
\end{itemize}
We  also  define   a  graded  $R$-linear  map 
\begin{eqnarray*}
dv:  ~~ \Lambda_n(V;R)\longrightarrow \Lambda_{n+1}(V;R),~~~ n\geq 0,
\end{eqnarray*}
by  letting  
\begin{eqnarray}\label{eq-2.2.1}
d  v(v_0v_1\ldots v_n)=
\sum_{i=0}^{n+1}  (-1)^i v_0\ldots v_{i-1} v v_i\ldots v_n.  
\end{eqnarray}
We extend (\ref{eq-1.1})  and   (\ref{eq-2.2.1})   linearly  over $R$.   
For any $u,v\in V$,  it  is  direct to  verify  
\begin{eqnarray}\label{eq-1.0}
\frac{\partial}{\partial u}\circ\frac{\partial}{\partial v}=-\frac{\partial}{\partial v}\circ\frac{\partial}{\partial u},~~~~~~  
d  u   \circ   dv  =  -   dv  \circ  du.   
\end{eqnarray}
In particular,   for  any  $v\in  V$,  
\begin{eqnarray}\label{eq-1.0111}
\frac{\partial}{\partial v}\circ \frac{\partial}{\partial v}=0,~~~~~~   dv \circ  dv  =0.   
\end{eqnarray}

\begin{proposition}
$\Lambda_*(V;R)$  is  
\begin{enumerate}[(i).]
\item
  a  $R_*(V)$-parametrized  DGA;      
  \item
     a   $R^*(V)$-parametrized  DGA.   
     \end{enumerate}
\end{proposition}

\begin{proof}
(i).  Condition~(1)  in  Definition~\ref{def-881}  follows  from  (\ref{eq-1.0111}).   Condition~(2)  in  Definition~\ref{def-881}  follows  from  (\ref{eq-1.0}).    Thus  $\Lambda_*(V;R)$  is  a  $R_*(V)$-parametrized  DGA  and   is  also  a  $R^*(V)$-parametrized  DGA.  
\end{proof}

Note  that   the   graded  free  $R$-module  $\Lambda_*(G;R)$  generated  by  the  elementary   paths  on   $G$
  is  a  graded   sub-$R$-module  of  $\Lambda_*(V;R)$.     
 By   Section~\ref{s1113},   
 for any $\alpha\in\wedge^p  R_*(V)$, where $p$  is odd,  any $\beta\in\wedge^r   R_*(V)$, where $r$  is even,   and       any $q=0,1,\ldots,p-1$,
   there     are  chain  complexes  
   \begin{eqnarray*}
  {\rm  Inf}_\bullet (G,\alpha,q)&:=& {\rm  Inf}_\bullet(\Lambda_*(G;R),\Lambda_*(V;R),\alpha,q),\\
    {\rm  Sup}_\bullet (G,\alpha,q)&:=& {\rm  Sup}_\bullet(\Lambda_*(G;R),\Lambda_*(V;R),\alpha,q).   
   \end{eqnarray*}
We  define    the {\it  $(\alpha,q)$-parametrized   path  homology groups}   of   $G$  by
\begin{eqnarray}\label{def-aa23}
H_*(G,\alpha,q):=H_*(\Lambda_*(G;R),\Lambda_*(V;R),\alpha,q).  
\end{eqnarray}
 We  have an  induced  morphism of  homology  groups
 \begin{eqnarray} \label{eq-x125}
  \beta_*: H_*(G, \alpha,q)\longrightarrow  H_*(G,\alpha,q-r).
 \end{eqnarray}
Similarly,  for any $\bar\alpha\in\wedge^p  R^*(V)$, where $p$  is odd,  any $\bar\beta\in\wedge^r   R^*(V)$, where $r$  is even,   and       any $q=0,1,\ldots,p-1$,
   there  is  are  chain  complexes  
   \begin{eqnarray*}
  {\rm  Inf}_\bullet (G,\bar\alpha,q)&:=& {\rm  Inf}_\bullet(\Lambda_*(G;R),\Lambda_*(V;R),\bar\alpha,q),\\
    {\rm  Sup}_\bullet (G,\bar\alpha,q)&:=& {\rm  Sup}_\bullet(\Lambda_*(G;R),\Lambda_*(V;R),\bar\alpha,q).   
   \end{eqnarray*}
We  define    the {\it  $(\bar\alpha,q)$-parametrized   path  up-homology groups}   of   $G$  by
\begin{eqnarray}\label{def-aa26}
H^*(G,\bar\alpha,q):=H^*(\Lambda_*(G;R),\Lambda_*(V;R),\bar\alpha,q).  
\end{eqnarray}
 We  have an  induced  morphism of  homology  groups
 \begin{eqnarray} \label{eq-x128}
  \bar\beta_*: H^*(G, \bar\alpha,q)\longrightarrow  H^*(G,\bar\alpha,q+r).
 \end{eqnarray}
 
 \begin{remark}
 If  $G$   is  transitive,  then  for any $\alpha\in\wedge^p   R_*(V)$, where $p$  is odd,  we  have  
 \begin{eqnarray*}
 {\rm  Inf}_\bullet (  G,\alpha,q)  =    {\rm  Sup}_\bullet (   G,\alpha,q) 
 \end{eqnarray*}
 and for any $\bar\alpha\in\wedge^p  R^*(V)$,  we  have  
  \begin{eqnarray*}
 {\rm  Inf}_\bullet (  G,\bar\alpha,q)  =    {\rm  Sup}_\bullet (   G,\bar\alpha,q).   
 \end{eqnarray*}
 \end{remark}
 
\begin{remark}
Even  for  a   digraph  with   a  small number of  vertices,     the  rank  of  the parametrized homology could  be  large.   We  do not  give  any  examples  explicitly  since  for  the  reason  of   computational  complexity.  
\end{remark}
 

\subsection{The  functoriality of  the  parametrized homology  of  digraphs}
 
 \begin{theorem}[Main result  I]
 \label{th-88978}
 Let  $G$   be  a  digraph.  Then  with respect to  morphisms  of  digraphs,   we  have 
 \begin{enumerate}[(i).]
 \item
 for any $\alpha\in\wedge^p   R_*(V)$, where $p$  is odd, any $\beta\in\wedge^r  R_*(V)$, where $r$  is even,   and any $q=0,1,\ldots,p-1$,
 both
 the $(\alpha,q)$-parametrized    homology  $H_*(G, \alpha,q)$   and
 the morphism $\beta_*$      in  (\ref{eq-x125})  are    functorial; 
\item
  for any $\bar\alpha\in\wedge^p   R^*(V)$, where $p$  is odd, any $\bar\beta\in\wedge^r  R^*(V)$, where $r$  is even,   and any $q=0,1,\ldots,p-1$,
 both
 the $(\alpha,q)$-parametrized    up-homology  $H^*(G, \alpha,q)$   and
 the morphism $\bar\beta_*$          in  (\ref{eq-x128})   are    functorial.  
 \end{enumerate}
 \end{theorem}
 
 \begin{proof}
 Let  $p:  G\longrightarrow   G'$   be  a  morphism  of  digraphs.   Then  $p_0$  induces  a  homomorphism  of  differential graded algebras  $\Lambda (p_0):  \Lambda_*(V;R)\longrightarrow  \Lambda_*(V;R)$.  
 The  restriction of   $\Lambda (p_0)$  induces a  well-defined  homomorphism  of  graded  $R$-modules  
 $\Lambda (p_0):    \Lambda_*(G;R)\longrightarrow  \Lambda_*(G';R)$.  
  This    induces   two  commutative  diagrams  
  \begin{eqnarray*}
  \xymatrix{
   {\rm  Inf}_\bullet (G,\alpha,q) \ar[r]^-{  \Lambda (p_0)} \ar[d]_-{\iota_G} & {\rm  Inf}_\bullet (G',\alpha,q)  \ar[d]^-{\iota_{G'}}\\
      {\rm  Sup}_\bullet (G,\alpha,q) \ar[r]^-{  \Lambda (p_0)}  & {\rm  Sup}_\bullet (G',\alpha,q)  
  }
  \end{eqnarray*}
  and 
   \begin{eqnarray*}
  \xymatrix{
   {\rm  Inf}_\bullet (G,\bar\alpha,q) \ar[r]^-{  \Lambda (p_0)} \ar[d]_-{\iota_G} & {\rm  Inf}_\bullet (G',\bar\alpha,q)  \ar[d]^-{\iota_{G'}}\\
      {\rm  Sup}_\bullet (G,\bar\alpha,q) \ar[r]^-{  \Lambda (p_0)}  & {\rm  Sup}_\bullet (G',\bar\alpha,q)  
  }
  \end{eqnarray*}
  where  the  horizontal  maps  are  chain  maps  and  the  vertical  maps  are  canonical  inclusions  of   chain  complexes.   By  Lemma~\ref{th-1.288},  we  obtain the  theorem.  
   \end{proof}

 The  next  corollary  follows   from Theorem~\ref{th-88978}.   

 \begin{corollary}\label{co-8912}
 Let  $G$  be  a  digraph.    Let   $\delta  G$  be  a  transitive-interior  of   $G$.   Then 
 \begin{enumerate}[(i).]
 \item
 for any $\alpha\in\wedge^p   R_*(V)$, where $p$  is odd, any $\beta\in\wedge^r  R_*(V)$, where $r$  is even,   and any $q=0,1,\ldots,p-1$,  we  have  a  functorial  commutative  diagram 
 \begin{eqnarray*}
 \xymatrix{
     H_*(\delta  G, \alpha,q)\ar[r]^{\beta_*}\ar[d]  & H_*(\delta  G,\alpha,q-r)\ar[d] \\
    H_*(G, \alpha,q)\ar[r]^{\beta_*}\ar[d]  & H_*(G,\alpha,q-r) \ar[d]\\
    H_*(\Delta  G, \alpha,q)\ar[r]^{\beta_*}   & H_*(\Delta  G,\alpha,q-r);    
 }
 \end{eqnarray*}
 \item
   for any $\bar\alpha\in\wedge^p   R^*(V)$, where $p$  is odd, any $\bar\beta\in\wedge^r  R^*(V)$, where $r$  is even,   and any $q=0,1,\ldots,p-1$,   we   have   a   functorial  commutative  diagram 
    \begin{eqnarray*}
 \xymatrix{
     H^*(\delta  G, \bar\alpha,q)\ar[r]^{\bar\beta_*}\ar[d]  & H^*(\delta  G,\bar\alpha,q+r)\ar[d] \\
    H^*(G, \bar\alpha,q)\ar[r]^{\bar\beta_*}\ar[d]  & H^*(G,\bar\alpha,q+r) \ar[d]\\
    H^*(\Delta  G, \bar\alpha,q)\ar[r]^{\bar\beta_*}   & H^*(\Delta  G,\bar\alpha,q+r).      
 }
 \end{eqnarray*}
 \end{enumerate}
 Here  all the  vertical  maps  are  induced  by  the  canonical  inclusions   $\delta  G\subseteq  G\subseteq  \Delta  G$.  
 \end{corollary}
 
 \begin{proof}
Consider  the   graded  free  $R$-modules   generated  by  the  elementary  $n$-paths  on  $\delta  G$,   $G$,   
 and   $\Delta   G$  respectively.   As   graded sub-$R$-modules,    we   have  
\begin{eqnarray*}
\Lambda_*(\delta  G;R)\subseteq  \Lambda_*(G;R) \subseteq   \Lambda_*(\Delta  G;R)  \subseteq  \Lambda_*(V;R).  
\end{eqnarray*}
  The  commutative   diagrams   in  (i)  and  (ii) follow   from  the  functorialities  in   Theorem~\ref{th-88978}~(i)  and  (ii)  respectively.    The  functorialities  in  (i)  and   (ii)  also   follow  from    Theorem~\ref{th-88978}~(i)  and  (ii)  respectively. 
 \end{proof}

 \begin{remark}
 In  Corollary~\ref{co-8912}~(i),   we  have  
 \begin{eqnarray*}
  {\rm  Inf}_\bullet (\delta  G,\alpha,q) &=&  {\rm  Sup}_\bullet (\delta  G,\alpha,q),  \\
     {\rm  Inf}_\bullet (\Delta  G,\alpha,q) &=&  {\rm  Sup}_\bullet (\Delta  G,\alpha,q).    
     \end{eqnarray*}
      In  Corollary~\ref{co-8912}~(ii),   we  have  
 \begin{eqnarray*}
  {\rm  Inf}_\bullet (\delta  G,\bar\alpha,q) &=&  {\rm  Sup}_\bullet (\delta  G,\bar\alpha,q),  \\
     {\rm  Inf}_\bullet (\Delta  G,\bar\alpha,q) &=&  {\rm  Sup}_\bullet (\Delta  G,\bar\alpha,q).    
     \end{eqnarray*}

 \end{remark}
 
 \subsection{Some  K\"unneth-type  formulae  for  the  parametrized homology  of  digraphs}

 Let  $V$  and  $V'$  be  two  disjoint  vertex  sets.    
Let  $G=(V,D,s,t)$  and   $G'=(V',D',s',t')$   be  two  digraphs.   The  {\it  join}  of   $G$  and   $G'$  is   a  digraph 
 \begin{eqnarray*}
 G*G'= (V\sqcup  V',  D*D',s*s', t*t')
 \end{eqnarray*}
 where        
 \begin{eqnarray*}
 D*D'=  D\sqcup  D'  \sqcup  \{(v,v')\mid   v\in D,  v'\in D'\}
 \end{eqnarray*}
 and    for  any  $e\in  D*D'$, 
 \begin{eqnarray*}
 (s*s')(e)&=&\begin{cases}
 s(e)  &  {\rm~if~}  e\in  D,\\
  s'(e)  &  {\rm~if~}  e\in  D',\\
  v  &{\rm~if~} e=(v,v')  {\rm~where~}v\in  V  {\rm~and~}v'\in  V',
 \end{cases}
 \\
  (t*t')(e)&=&\begin{cases}
 t(e)  &  {\rm~if~}  e\in  D,\\
  t'(e)  &  {\rm~if~}  e\in  D',\\
  v'  &{\rm~if~} e=(v,v')  {\rm~where~}v\in  V  {\rm~and~}v'\in  V'.  
 \end{cases}
 \end{eqnarray*} 
 Let  $n$  and  $m$  be  nonnegative  integers.  We  have  an  $R$-bilinear  map
 \begin{eqnarray}\label{eq-qwe889}
 *:  \Lambda_n(V;R)\times  \Lambda_m(V';R)\longrightarrow  \Lambda_{n+m+1}(V\sqcup  V';R)
 \end{eqnarray}
 sending  a  pair $(v_0\ldots  v_n,  v'_0\ldots  v'_m)$   of  an  elementary  $n$-path  $v_0\ldots  v_n$  on  $V$  and  an  elementary  
 $m$-path   $ v'_0\ldots  v'_m$   on  $V'$  to  the   elementary  $(n+m+1)$-path  $v_0\ldots  v_n  v'_0\ldots  v'_m$   
 on  $V\sqcup  V'$.    The restriction  of  (\ref{eq-qwe889})  to  the  paths on  $G$  and  $G'$  induces  an   $R$-bilinear  map 
  \begin{eqnarray}\label{eq-qwe8898}
 *:  \Lambda_n(G;R)\times  \Lambda_m(G';R)\longrightarrow  \Lambda_{n+m+1}(G*G';R). 
 \end{eqnarray}
 For  any  $x\in  R_*(V\sqcup  V')$,  there exist  a  unique  $y\in  R_*(V)$  and  a  unique $y'\in  R_*(V')$  such  that   $x=y+y'$.    For  any  $\xi\in  \Lambda(V;R)$  we  have 
$
 x(\xi)=y_1(\xi)
$
 and  for  any  $\xi'\in  \Lambda(V';R)$  we  have  $x(\xi)=y_2(\xi)$.  
 The  next  theorem  follows  from  Lemma~\ref{le-ceh2}.  
 
 \begin{theorem}[Main  result  II]
 \label{th-981328}
 Let  $G=(V,D,s,t)$  and   $G'=(V',D',s',t')$ be  two  digraphs  such  that   $V$  and   $V'$   are   disjoint.  Let  $n$  be  a  nonnegative  integer.   
 Then  for  any  $x\in  R_*(V\sqcup  V')$,    there  is  a  natural  short  exact  sequence 
  \begin{eqnarray}
&0\longrightarrow  \bigoplus_{t+s=n+1} H_t(G, x,0)\otimes H_s(G', x,0)\longrightarrow H_n(G*G',  x,0)
 \nonumber\\
 &\longrightarrow \bigoplus_{t+s=n+2 }{\rm Tor}_R(H_t(G, x,0),H_{s-1}(G', x,0))\longrightarrow 0.
 \label{eq-5.a9998aaaaa}
 \end{eqnarray}
 \end{theorem}
 
 \begin{proof}
 We  let $A=\Lambda_*(V;R)$  and  $B=\Lambda_*(G;R)$  in    Lemma~\ref{le-ceh2}.  For  any  nonnegative  integer   $n$,  it  is  direct  to see  that 
 \begin{eqnarray*}
 \Lambda_{n}(V\sqcup  V';R)=\bigoplus_{t+s=n+1\atop  t,s\geq 0}  \Lambda_t(V;R)\otimes  \Lambda_s(V';R)  
 \end{eqnarray*}
 whose  restriction  to  the  paths on  $G$  and  $G'$  gives 
  \begin{eqnarray*}
 \Lambda_{n}(G*   G';R)=\bigoplus_{t+s=n+1\atop  t,s\geq 0}  \Lambda_t(G;R)\otimes  \Lambda_s(G';R).    
 \end{eqnarray*}
The  natural  short  exact  sequence   (\ref{eq-5.a9998aaaaa})  follows  from Lemma~\ref{le-ceh2}.    
  \end{proof}

 The  next  corollary  follows  from  Theorem~\ref{th-981328}.  
 
 \begin{corollary}\label{co-uu76}
  Let  $G=(V,D,s,t)$  and   $G'=(V',D',s',t')$ be  two  digraphs  such  that   $V$  and   $V'$   are   disjoint.  Let  $n$  be  a  nonnegative  integer.    Let  $\delta  G$  be  a  transitive-interior  of   $G$  and  let  $\delta  G'$  
  be  a  transitive  interior  of   $G'$.  
 Then  for  any  $x\in  R_*(V\sqcup  V')$,    there  is  a     commutative  diagram  
 \begin{eqnarray*}
 \xymatrix{
 0\ar[r]  & \bigoplus_{t+s=n+1} H_t(\delta  G, x,0)\otimes H_s(\delta  G', x,0)\ar[r]  \ar[d]   &H_n(\delta  G*  \delta  G',  x,0) \ar[r] \ar[d]& \\
0\ar[r]  & \bigoplus_{t+s=n+1} H_t(G, x,0)\otimes H_s(G', x,0)\ar[r]  \ar[d]   &H_n(G*G',  x,0) \ar[r] \ar[d]& \\
 0\ar[r]  & \bigoplus_{t+s=n+1} H_t(\Delta  G, x,0)\otimes H_s(\Delta  G', x,0)\ar[r]     &H_n(\Delta  G*  \Delta  G',  x,0)\ar[r]  &
 }\\
  \xymatrix{
 \ar[r]   & \bigoplus_{t+s=n+2 }{\rm Tor}_R(H_t(G, x,0),H_{s-1}(G', x,0))\ar[r]  \ar[d]  & 0\\
 \ar[r]  & \bigoplus_{t+s=n+2 }{\rm Tor}_R(H_t(G, x,0),H_{s-1}(G', x,0))\ar[r]  \ar[d]  & 0\\
 \ar[r]   & \bigoplus_{t+s=n+2 }{\rm Tor}_R(H_t(G, x,0),H_{s-1}(G', x,0))\ar[r]    & 0
 }
 \end{eqnarray*}
 such that  each  row  is  a  natural  short  exact  sequence.  
 \end{corollary}
 
 \begin{proof}
 The  commutative  diagram  follows  by applying  the  inclusions  $\delta  G\subseteq  G\subseteq  \Delta   G$ to   the  naturality of  the   short  exact sequence  (\ref{eq-5.a9998aaaaa}).  
 \end{proof}

  \bigskip

Shiquan Ren (first author)

Address:  School of  Mathematics and Statistics,  Henan  University,   Kaifeng,  Henan,  475004,  P. R. China

E-mail:  renshiquan@henu.edu.cn

\bigskip

Chong Wang   (corresponding author)

Address:    School
of Mathematics and Statistics, Cangzhou Normal University, Cangzhou, Hebei, 061000, P. R. China

E-mail:  wangchong\_618@163.com


\begin{thebibliography}{99}







\bibitem{h1}
S. Bressan, J. Li, S. Ren   and  J.  Wu,  \emph{The  embedded homology of hypergraphs and applications},     Asian J.    Math. {\bf 23}(3)   (2019),  479-500.


 


\bibitem{lin1}
A.   Grigor'yan,  Y.   Lin, S.-T. Yau,  \emph{Torsion of digraphs and path complexes},  arXiv: 2012.07302v1,  2020.

\bibitem{lin2}
A. Grigor'yan,  Y.  Lin, Y. Muranov,  S.-T. Yau,  \emph{Homologies of path complexes and digraphs},  arXiv: 1207.2834,  2013.

\bibitem{lin3}
A. Grigor'yan,  Y.  Lin, Y.  Muranov,  S.-T. Yau,  \emph{Homotopy  theory for  digraphs},   Pure and Applied Mathematics Quarterly   {\bf 10}(4),   619-674,  2014.


\bibitem{lin4}
A.  Grigor'yan,  Y.  Lin, Y.  Muranov,  S.-T. Yau,  \emph{Cohomology of digraphs and (undirected) graph},  Asian J. of  Math.  {\bf 15} (5),   887-932,  2015.



\bibitem{lin6}
A.  Grigor'yan,   Y. Lin,  Y. Muranov,  S.-T. Yau,  \emph{Path complexes and their homologies},   Journal of Mathematical Sciences{\bf 248} (5),   564-599,  2020.


\bibitem{lin5}
A.  Grigor'yan,   Y.  Muranov,  S.-T. Yau,  \emph{Homologies of digraphs and K\"{u}nneth formulas},   Communications in Analysis and Geometry  {\bf 25},   969-1018,  2017.












\bibitem{hatcher}
A.   Hatcher,  Algebraic topology. Cambridge University Press,  2002.




\bibitem{wrl}
C. Wang,  S. Ren,   J. Liu,  \emph{A  K\"unneth formula  for finite sets},   Chin.  Ann.  of  Math.  Ser. B {\bf 42}(6), 801-812, 2021.



  \end{thebibliography}
\end{document}